\documentclass{elsart3-1}

\usepackage{amssymb}

\usepackage[english,francais]{babel}

\newtheorem{theorem}{Theorem}[section]

\newtheorem{lemma}[theorem]{Lemma}

\newtheorem{e-proposition}[theorem]{Proposition}

\newtheorem{e-definition}[theorem]{Definition\rm}

\newtheorem{remark}{\it Remark\/}

\renewcommand{\theequation}{\arabic{equation}}

\def\og{\leavevmode\raise.3ex\hbox{$\scriptscriptstyle\langle\!\langle$~}}

\def\fg{\leavevmode\raise.3ex\hbox{~$\!\scriptscriptstyle\,\rangle\!\rangle$}}

\begin{document}

\begin{frontmatter}













\selectlanguage{english}

\title{Stochastic Loewner evolution in multiply connected domains}






\selectlanguage{english}
\author[authorlabel1]{Robert~O. Bauer\thanksref{1}},
\ead{rbauer@math.uiuc.edu}
\author[authorlabel2]{Roland~M. Friedrich\thanksref{2}}
\ead{rolandf@ias.edu}

\thanks[1]{Research supported by NSA grant H98230-04-1-0039}
\thanks[2]{Research supported by NSF grant DMS-0111298}

\address[authorlabel1]{Department of Mathematics, University of
Illinois at Urbana-Champaign, 1409 West Green Street, Urbana, IL 61801, USA}
\address[authorlabel2]{Institute for Advanced Study, Princeton,
NJ 08540, USA}

\begin{abstract}


We construct radial stochastic Loewner evolution in multiply connected domains, choosing the unit disk with concentric circular slits as a family of standard domains. The natural driving function or input is a diffusion on the associated moduli space. The diffusion stops when it reaches the boundary of the moduli space. We show that for this driving function the    family of random growing compacts has a phase transition for $\kappa=4$ and $\kappa=8$, and that it satisfies locality for $\kappa=6$.

\end{abstract}

\end{frontmatter}

\selectlanguage{english}


\section{Introduction}

\label{}

The stochastic Loewner evolution was introduced by Schramm in \cite{schramm:2000}. The Loewner equation shows how to encode a slit in the unit disk which starts at the boundary circumference in terms of a motion on the boundary. For a random slit such as an interface created by a statistical mechanics model at criticality, Schramm derived just  from a few properties anticipated for the scaling limit of such a model---conformal invariance, a Markovian-type property, symmetry, continuity---that the random motion on the boundary encoding the random slit has to be a linear time-change of a standard Brownian motion. This new family of processes, $SLE_{\kappa}$, $(\kappa>0)$, is studied in \cite{lsw:2001a}, \cite{lsw:2001b}, \cite{lsw:2003}, \cite{rohde.schramm:2001}, and references therein. In particular, it was shown that a phase transition occurs at $\kappa=4$ where the slit ceases to be a simple curve and at $\kappa=8$, where the slit becomes space-filling. For certain values of $\kappa$, $SLE_{\kappa}$ has special properties---locality for $\kappa=6$, and the restriction property for $\kappa=8/3$.

Since the statistical mechanics models whose scaling limit is (or is conjectured to be) described by $SLE_{\kappa}$ are readily defined for non-simply connected domains the question arises if the stochastic Loewner evolution can be extended to the non-simply connected case. The Loewner equation for simply connected domains has an analogue for multiply connected  domains which we call the Komatu-Loewner equation. We derive a version of a result by Komatu for a different class of standard domains, \cite{komatu:1950}.  Cutting a non-simply connected domain along a Jordan arc changes the conformal type of the domain. Because of this, the Kommatu-Loewner equation is not enough to reverse the process, i.e. to begin with a driving function on the boundary and solve the equation---the coefficients in the equation depend on the conformal type of the domain. However, using Schiffer variations we obtain an associated system of equations for the conformal invariants. These equations, together with the Komatu-Loewner equation, allow us to reverse the process.

In the case of a random slit, the expected  properties of the  slit Schramm used for simply connected domains only allow us to conclude that the Komatu-Loewner equation should be driven by a diffusion on the space of conformal invariants, the moduli space. We define a class of diffusions on the moduli space for which the appropriate phase transitions occur at $\kappa=4$ and $\kappa=8$, and which satisfies locality at $\kappa=6$.

\section{Komatu-Loewner equation}

Denote by $D$ a domain in the complex plane bounded by a finite number $n\ge2$ of
proper continua. Then $D$ is conformally equivalent to the unit disk cut along $n-1$ disjoint concentric circular slits. We call such a domain
a standard domain. If $D$ is a standard domain we denote the boundary
components by
\begin{equation}
  C^{(j)}: |z|=m_j,\ \theta_j\le \arg z\le \theta_j',\ 1\le j\le n-1 ,\quad C^{(n)}: |z|=1.
\end{equation}
Let $\tilde{D}$ be a domain in the $w$-plane
obtained from a standard domain by cutting
along a slit (Jordan arc) $\Gamma$  which starts from a point on the
exterior  boundary component $|w|=1$ and avoids $0$. Denote $w=f(z)$
the unique conformal map from a standard domain
$D$ onto $\tilde{D}$ such that $f(0)=0$, $f'(0)>0$, and such that the peripheral boundary circumferences correspond to each
other. Let $T=-\ln f'(0)>0$. If we delete from the boundary of $\tilde{D}$ a
sub-arc of $\Gamma$ with interior endpoint common with it, then there is
a unique conformal map
\begin{equation}
  w=h(w_t,t),\quad h(0,t)=0, \quad h'(0,t)=e^{-t}
\end{equation}
which maps a standard domain onto the domain thus obtained. Here $t\le T$ and the peripheral
boundary components are to correspond to each other. Define $D_t$ by $h(D_t,t)=\tilde{D}$.
Then $D_t$ is a domain of the same type as $\tilde{D}$, i.e. obtained from a
standard domain by cutting along a slit $\Gamma_t$. By the monotony property
of the derivative $h'(0,t)$ the points on $\Gamma$ are in one-to-one
correspondence with the values of the parameter $t$, ranging over the interval
$0\le t\le T$. Let the conformal map from $D$ onto $D_t$ be denoted by
\begin{equation}
  w_t=f(z,t),\quad f(0,t)=0, \quad f'(0,t)=e^{t-T}.
\end{equation}
Then $
  f(z)=f(z,0)=h(f(z,t),t)$,
and also $f(z,T)=z$. Further, let the boundary components of the circular slit disk
$D_t+\Gamma_t$ be denoted by
\begin{equation}
  C^{(j)}_t: |w_t|=m_j,\ \theta_j(t)\le \arg w_t\le \theta_j'(t),\ 1\le j\le
  n-1, \quad C^{(n)}_t: |w_t|=1,
\end{equation}
and the starting point on $C^{(n)}_t$ of the slit $\Gamma_t$ be
$\overline{\gamma}(t)=e^{-i\theta(t)}$. Denote $G(u,w_t;t)$ the Green function
of
the circular slit disk $D_t+\Gamma_t$ in the $u$-plane with pole at $w_t$, and denote a harmonic conjugate to $G(u,w_t;t)$ with respect to $w_t=x_t+i y_t$ by
\begin{equation}
        H(u,w_t;t)=\int_0^{w_t}\left(\frac{\partial G(u,w_t;t)}{\partial x_t}\ dy_t-\frac{\partial G(u,w_t;t)}{\partial y_t}\ dx_t\right),
\end{equation}
and set $
        F(u,w_t;t)=G(u,w_t;t)+i H(u,w_t;t)$.
Denote
$\omega_j(w_t;t)$ the harmonic measure of $C^{(j)}_t$ at $w_t$ with respect to
$D_t+\Gamma_t$ and let
\begin{equation}
        R_j(w_t;t)=\frac{1}{2\pi}\int_{C_t^{(j)}}\frac{\partial F(u,w_t;t)}{\partial\nu}\ ds,
\end{equation}
where $\partial/\partial \nu$ denotes differentiation at the
boundary point $u$ in the direction of the inner normal. We have $\Re(R_j)=\omega_j$ and we may assume that $\Im(R_j(0;t))=0$ for $j=1,\dots,n-1$.

\begin{theorem}[Radial Komatu-Loewner equation]
  The family $\{f(z,t), 0\le t\le T\}$ satisfies the equation
  \begin{equation}\label{E:mult-loewner}
    \frac{\partial\ln f(z,t)}{\partial t}
    =\frac{\partial F(\overline{\gamma}(t),w_t;t)}{\partial\nu}
    +\sum_{j=1}^{n-1} R_j(w_t;t)\frac{d\ln m_j(t)}{dt},
  \end{equation}
  with initial condition $f(z,T)=z$ for all $z\in D\backslash\{0\}$.
\end{theorem}

Note that, while the first and also the second term on the right in equation (\ref{E:mult-loewner}) are multiple-valued, their sum is single-valued.
In the simply connected case the sum on the right-hand-side of equation
    (\ref{E:mult-loewner}) disappears and the radial Komatu-Loewner equation reduces to the usual radial Loewner equation,
    \begin{equation}
      \frac{\partial\ln f(z,t)}{\partial t}=-
      \frac{f(z,t)+\overline{\gamma}(t)}{f(z,t)-\overline{\gamma}(t)}.
    \end{equation}

\begin{remark} An alternative class of standard domains is provided by
    circular slit annuli. For that case Komatu \cite{komatu:1950} derived an
    equation analogous to  (\ref{E:mult-loewner}). In the case $n=2$ it can be
    written in terms of elliptic functions, \cite{komatu:1943}.
\end{remark}

If
\begin{equation}
  {\bf P}={\bf P}_t=[p_{j,k}(t)]_{j,k=1}^{n-1},\quad
  p_{j,k}(t)=\int_{C_t^{(k)}}\frac{\partial \omega_j(w_t;t)}{\partial\nu}\ ds
\end{equation}
is the period matrix of $D_t$ and $ \lambda_j=-\ln m_j$,
${\bf\lambda}^T=[\lambda_1,\dots,\lambda_{n-1}]$, then
\begin{equation}
  \lambda^T=2\pi\omega^T(0){\bf P}^{-1},
\end{equation}
where $\omega^T=(\omega_1,\dots,\omega_{n-1})$ is the vector of harmonic measures.
Under a Schiffer variation at the boundary point $\overline{\gamma}(t)$
\begin{equation}
  \frac{\partial}{\partial
  \epsilon}|_{\epsilon=0}2\pi\omega^T(0){\bf P}^{-1}=2\pi
  \frac{\partial\omega(\overline{\gamma}(t),t)}{\partial\nu}^T{\bf P}^{-1},
\end{equation}
see \cite{schiffer:1946}. Furthermore, the domain constant $d_n(t)$ defined by
\begin{equation}
  \Re{(\ln f^{-1}(w_t,t))}=\ln|w_t|+d_n(t)+O(|w_t|)
\end{equation}
is given by $d_n(t)=T-t$. Under a Schiffer variation we have
$
  \frac{\partial}{\partial
  \epsilon}|_{\epsilon=0}d_n(t)=-1$, see \cite{schiffer:1946}.
We find
\begin{equation}
  \sum_{j=1}^{n-1} R_j(w_t;t)\frac{d\ln m_j(t)}{dt}= 2\pi
  \frac{\partial\omega(\overline{\gamma}(t),t)}{\partial\nu}^T{\bf
  P}_t^{-1}{\bf R}(w_t;t).
\end{equation}
If $n>1$, then the moduli space of an $n$-connected domain with one interior marked point is $3n-4$ dimensional. Denote $M_l(t),1\le l\le
3n-4$ the moduli of the domain $D_t+\Gamma_t$ with the origin as marked point.

\begin{lemma}
  There exist differentiable vector-fields $Y_l=Y_l(u,M_1,\dots, M_{3n-4})$, $1\le l\le 3n-4$,
  where $u$ is complex of norm 1 and the $M_k$ real such that if
  $M_1(t),\dots,M_{3n-4}(t)$ solve
  \begin{equation}\label{E:moduli-motion}
    \dot{M}_l(t)=Y_l(\overline{\gamma}(t),M_1(t),\dots,M_{3n-4}(t)), 1\le l\le
    3n-4,
  \end{equation}
with initial values the parameters of the domain $D_0+\Gamma_0$,
  then $\{f(z,t),0\le t\le T\}$ solves (\ref{E:mult-loewner}) with the Green
  function and harmonic measure at time $t$ calculated from the parameters
  $M_k(t)$.
\end{lemma}

We can now reverse the procedure. Instead of beginning with a Jordan arc
$\Gamma$ that leads to a continuous curve $t\mapsto\overline{\gamma}(t)$ on
the boundary circumference we begin with the curve $\overline{\gamma}(t)$.

\begin{theorem}
  For a standard domain $D$ with parameters $M_1,\dots, M_{3n-4}$ and a continuous function $t\mapsto\overline{\gamma}(t)$ with values in
  the unit circle we can solve (\ref{E:moduli-motion}) with initial values $M_1,\dots, M_{3n-4}$. The solution exists on $(0,T]$. If $T<\infty$, then there is a $j$ such that $\lim_{t\nearrow T} m_j(t)=1$. For $z\in D$ and $t<T$, let $g_t(z)$ be the solution
of (\ref{E:mult-loewner}) with $g_0(z)=z$. The solution exists up to  a time $T_z\in(0,T]$. If $T_z<T$, then $\lim_{t\nearrow T_z}(g_t(z)-\overline{\gamma}(t))=0$. If $K_t$ denotes the closure of the set of points $z$ for which $T_z\le t$, then $g_t$ maps $D\backslash K_t$ conformally onto a standard domain $D_t$.
\end{theorem}


\section{Diffusion on moduli space}

To define $SLE_{\kappa}$ on $D$ we first define a diffusion on an appropriate moduli space. This is the moduli space ${\mathfrak M}_n$ of an $n$-connected domain with one interior point and one boundary point marked. It has dimension $3n-3$.

Let $\kappa>0$ and denote $B_t$ a standard one dimensional Brownian motion. Define $l(e^{i\theta};t)$ by
\begin{equation}
        l(e^{i\theta};t)=\lim_{\varphi\to \theta}\left(\frac{\partial F(e^{i\theta},e^{i\varphi};t)}{\partial\nu}-\frac{e^{i\theta}+e^{i\varphi}}{e^{i\theta}-e^{i\varphi}}\right).
\end{equation}
Consider the system
\begin{eqnarray}\label{E:diffusion}
        d\theta(t)&=&\frac{1}{i}\left(l(\overline{\gamma}(t);t))\ dt
    +\sum_{j=1}^{n-1} R_j(\overline{\gamma}(t);t)\ d\ln m_j(t)\right)+\sqrt{\kappa}dB_t\\
        dM_l(t)&=&Y_l(\overline{\gamma}(t), M_1(t),\dots, M_{3n-4}(t)),\quad l=1,\dots,3n-4,
\end{eqnarray}
where $\theta(0)=0$, $\overline{\gamma}(t)=e^{i\theta(t)}$, and $M_l(0)=M_l$. We note again that the combination of terms in brackets in equation (\ref{E:diffusion}) is single-valued. The solution exists up to a stopping time $\tau$. On the set $\tau<\infty$ we have for some $j$, $\lim_{t\nearrow\tau}m_j(t)=1$. We call a solution to this system a {\em Schiffer diffusion (on moduli space)\ } and the random family of growing compacts $K_t$ associated to a Schiffer diffusion via (\ref{E:mult-loewner}) {\em radial $n-SLE_{\kappa}$}. Alternatively, we may call the family of random conformal maps $g_t$ $n-SLE_{\kappa}$. The following results show that the Schiffer diffusion has key properties we expect from our knowledge of SLE in simply connected domains.

\begin{theorem}
        If $\kappa\le 4$, then $K_t$ is  a.s. a simple curve for $t<\tau$.      If $\kappa>4$  then $K_t$ is a.s. not simple.
\end{theorem}

This follows by considering the motion of points on the boundary.

\begin{remark}
        We do not know if $K_t$ can hit a circular slit in finite time with positive probability.
\end{remark}

\begin{theorem}
        If $\kappa=6$, then $n-SLE_{\kappa}$ satisfies locality. 
\end{theorem}

The proof is similar to the simply connected case. Our definition of the Schiffer diffusion was made with the previous two theorems in mind. We have not been able to confirm the restriction property for $\kappa=8/3$. Instead we only obtained a family of martingales related to the restriction property. We will investigate in a forthcoming paper if the Schiffer diffusion is the essentially only class of diffusions for which these properties hold, and also if locality holds for $\kappa=8/3$.

In a multiply connected domain there are two flavors of chordal SLE: Growing slits from a point on one boundary component to another point on the same boundary component, or growing slits from a point on one boundary component to a point on another boundary component. For the former we can repeat the above construction using as standard domain the upper half-plane with a finite number of horizontal slits, the compacts $K_t$ growing toward $\infty$. For the latter the construction can be based on Komatu's equation \cite{komatu:1950}.







\section*{Acknowledgements}


The authors would like to acknowledge the hospitality of the Institut f\"ur Stochastik at the University of Freiburg and of the Isaac Newton Institute where part of this work was done.

\end{document}